\documentclass{article}
\usepackage[utf8]{inputenc}
\usepackage{amssymb,amsmath,latexsym}
\usepackage[english]{babel}

\title{On Certain Reciprocal sums}
\author{Soumyadip Sahu\\soumyadip.sahu00@gmail.com}
\date{}

\begin{document}
\maketitle
\begin{abstract}
In this note we associate a sequence of non-negative integers to any convergent series of positive real numbers and study this sequence for the series $\sum_{n \geq 1} n^{-k}$ where $k$ is an integer $\geq 2$. 
\end{abstract}
\section{Introduction}
Let $(x_i)_{i=1}^\infty$ be a sequence of positive real numbers such that $\sum_{i=1}^\infty x_i$ converges. Given such a sequence one can associate a sequence of non-negative integers $(a_n)_{n=1}^{\infty}$ by defining \[a_n=\Big[\frac {1} {\sum_{i=n+1}^\infty x_i}\Big] \] where $[x] = \text{the largest integer} \leq x$, for a real number $x$. \\~\\
This problem has been studied for some special class of sequences. For example, Ohtsuka and Nakamura \cite{ohtsuka} derived a formula for $x_i=\frac {1} {F_i}$, where $F_i$ denotes the $i^{th}$ Fibonacci number. Since then several results have been discovered about the case in which $x_i$ s are reciprocals of a sequence given by linear recurrence relations (for example see \cite{others}).\\
In this note we consider the case $x_i=i^{-k}$ where $k$ is a positive integer $\geq 2$.\\
The first theorem that we prove is following:\\~\\
\textbf{Theorem 1.1 :} Let $k$ be an integer $\geq 2$. Then there is  a polynomial $f(X) \in \mathbb{Q}[X]$ of degree $(k-1)$ (unique upto the constant term) and an integer $N_0$ (depending on $f$) such that\\
\[f(n)  < \frac {1} {\sum_{i=n+1}^{\infty} i^{-k}} < f(n)+1\]
holds for all $n \geq N_0$.\\~\\
In section-2 we prove a lemma which is central to our treatment. In section-3 we prove theorem 1.1 and indicate how to compute a closed form formula for $a_n$  and compute it for $k=2,3,4,5$. In section-4 we prove a generalization which is as follows :\\~\\
Let $P(X)$ be a polynomial over $\mathbb{R}$  of degree $\geq 2$ such that the leading co-efficient is positive. Let $i_0 \in \mathbb{R}$ be large enough so that $P(x) > 0$ for all $x > i_0$. Put $x_i = \frac {1} {P(i+i_0)}$ for all $i \geq 1$. Clearly $\sum_{i \geq 0} x_i < \infty$. Then we have an analogus result :\\~\\
\textbf{Theorem 1.2 :} There is a polynomial $f(X) \in \mathbb{R}[X]$ depending on $P$ (unique upto constant term), an integer $N_0$ depending on $f$ and $i_0$ so that degree of $f$ is $(k-1)$ and 
\[f(n+i_0) < \frac {1} {\sum_{i=n+1}^{\infty} x_{i}} < f(n+i_0)+1\]
holds for all $n \geq N_0$.

\section{An Important Lemma}
We begin by proving a useful lemma.\\~\\
\textbf{Lemma 2.1:} Fix an integer $k$, $k\geq 2$. Let $x_0, x_1,\cdots,x_{k-1}$ be $k$ unknowns. Consider $F(X,x_0,\cdots,x_{k-1})=x_0X^{k-1}+x_1X^{k-2}+\cdots+x_{k-2}X+x_{k-1} \in \mathbb{R} [X,x_0,\cdots,x_{k-1}]$.
Let $g(X) \in \mathbb{R}[X]$ be a polynomial of degree $k$. Assume that \[g(X+1) = a_0X^k + \cdots + a_k.\]
Put
\[F((X+1),x_0,\cdots,x_{k-1}) 
=x_0X^{k-1}+(x_1+y_1)X^{k-2}+\cdots+(x_{k-2}+y_{k-2})X+(x_{k-1}+y_{k-1}),\]
\[H(X)=F((X+1),x_0,\cdots,x_{k-1})F(X,x_0,\cdots,x_{k-1})\\
      = p_0X^{2k-2}+\cdots+p_{2k-3}X+p_{2k-2},\]
\[G(X)= g(X+1)(F((X+1),x_0,\cdots,x_{k-1})- F(X,x_0,\cdots,x_{k-1}))=q_0X^{2k-2}+\cdots+q_{2k-3}X+q_{2k-2}\]
where $y_i, p_j, q_l \in \mathbb{R} [x_0,\cdots,x_{k-1}]$ for all $1 \leq i \leq k-1$, $0 \leq j \leq 2k-2$ and $0 \leq l \leq 2k-2$.\\
Consider the system of $k$ equations 
$$p_i=q_i ,\hspace{1mm}\forall \hspace{1mm}0 \leq i \leq k-1$$
 in $k$ unknowns $x_0,\cdots,x_{k-1}$.\\
This system of equations has a unique solution $(c_0,\cdots,c_{k-1}) \in \mathbb{R}^{k}$ with $c_0 \neq 0$.\\
Further, if $g$ is defined over $\mathbb{Q}$ then $(c_0, \cdots, c_{k-1}) \in \mathbb{Q}^{k}$.\\~\\
\textbf{Proof:} First notice that the coefficient of $X^k$ in $F((X+1),x_0,\cdots,x_{k-1})$ is indeed $x_0$ and degrees of $G(X), H(X)$ in $X$ are indeed at most $(2k-2)$. Hence the hypothesis of the lemma is justified.\\ 
Now an application of binomial theorem gives \[y_i=\binom {k-i} {1}x_{i-1}+\binom {k-i+1} {2}x_{i-2} + \cdots + \binom {k-1} {i}x_0 \tag{2.1}\] for each $1 \leq i \leq k-1$.\\
A direct calculation gives 
\[p_j=\sum_{r=0}^{j} x_r(x_{j-r}+y_{j-r}) \tag{2.2}\]
for all $0 \leq j \leq k-1$ where $y_0=0$.\\
Similarly, 
\[q_l=\sum_{r=0}^{l} a_{r}y_{l-r+1} \tag{2.3}\]
for all $0 \leq l \leq k-1$ where $y_k=0$.\\
Thus $p_0=x_0^2$ and $q_0=y_1=a_o\binom {k-1} {1} x_0=a_0(k-1)x_0$.\\
Hence $x_0=a_0(k-1)$ is a solution to $p_0=q_0$. Clearly $a_0(k-1) \neq 0$.\\
Now notice that $p_i$ depends only on $\{x_0,\cdots,x_i,y_0,\cdots,y_i\}$.\\
By $(2.1)$, this observation implies $p_i$ depends only on $\{x_0,\cdots,x_i\}$. This holds for all $0 \leq i \leq k-1$. Similarly $q_i$ depends only on $\{y_1,\cdots,y_{i+1}\}$ i.e. only on $\{x_0,\cdots,x_i\}$. This is true for all $0 \leq i \leq k-1$.\\ 
So we can take an inductive approach to solve the system of equations.\\
We have already found a $c_0$ (namely $a_0(k-1)$) such that $x_0=c_0$ solves $p_0=q_0$ and $c_0 \neq 0$. Note that this is the only non-zero solution to $p_0 = q_0$ and if $a_0 \in \mathbb{Q}$ then $c_0 \in \mathbb{Q}$.\\
Assume that we have found $(c_0,\cdots,c_i) \in \mathbb{R}^{i+1}$ such that $c_0 \neq 0$ and this is the unique tuple solving the system of equations
\[p_0=q_0,\cdots,p_i=q_i\] for some $i$ in the range $0 \leq i \leq k-2$. Further if $g$ is defined over $\mathbb{Q}$ then $\mathbb{Q}^{i+1}$. \\
Now the goal is to find a $c_{i+1} \in \mathbb{R}$ such that $(c_0,\cdots,c_{i+1})$ solves $p_{i+1}=q_{i+1}$.\\
We consider two cases:\\~\\
\textbf{Case I: } $i < k-2$\\
Here $i+1 \leq k-2$.\\
Coefficient of $x_{i+1}$ in $p_{i+1}$ is $=2x_0$ (one $x_0$ arises from term $x_0(x_{i+1}+y_{i+1})$ and other $x_0$ arises from the term $x_{i+1}(x_0+y_0)$).(Follows from $(2.1)$ and $(2.2)$)\\ 
Coefficient of $x_{i+1}$ in $q_{i+1}$ is = $a_0$(coefficient of $x_{i+1}$ in $y_{i+2})= a_0\binom {k-(i+2)} {1}= a_0(k-(i+2))$ (follows from $(2.1)$ and $(2.2)$).\\
Hence $p_{i+1}=q_{i+1}$ can be rewritten as  \[\{2x_0 - a_0(k - i -2)\}x_{i+1} = \text{some polynomial in}\; x_0,\cdots,x_i \tag {2.4}\]
Note that $\{2c_0 - a_0(k - i -2)\} = a_0(k+i) \neq 0$.
Thus we can put $x_0=c_0,\cdots,x_i=c_i$ in $(2.4)$ and solve for $x_{i+1}$ to get a tuple $(c_0, \cdots, c_{i+1}) \in \mathbb{R}^{i+2}$ which is a solution for the system of equations \[p_0 = q_0, \cdots, p_{i+1} = q_{i+1}.\]
Now if $(d_0,\cdots, d_{i+1})$ is another solution with $d_0 \neq 0$ then by induction hypothesis $(c_0, \cdots, c_{i}) = (d_0, \cdots, d_i)$. From $(2.4)$ it follows that $d_{i+1} = c_{i+1}$. Hence the uniqueness.\\
If $g$ is defined over $\mathbb{Q}$ then by induction hypothesis $(c_0, \cdots, c_i) \in \mathbb{Q}^{i+1}$. Since $p_{i+1}$ and $q_{i+1}$ are defined over $\mathbb{Q}$ using $(2.4)$ we conclude that $c_{i+1} \in \mathbb{Q}$.\\  
So for this case we are done.\\~\\
\textbf{Case II: } $i=k-2$\\
This is essentially similar to previous case. Only difference is coefficient of $x_{k-1}$ in $q_{k-1}$ is $0$.\\
So $p_{k-1}=q_{k-1}$ can be rewritten as $2x_0x_{k-1}=$ some polynomial in $x_0,\cdots,x_{k-2}$ and from here the arguments of previous case goes through since $c_0 \neq 0$.\\
Thus inductively we can find $(c_0,\cdots,c_{k-1}) \in \mathbb{R}^{k}$ such that this tuple is the unique solution to the system of equations under consideration with $c_0 \neq 0$. Further if $g$ is defined over $\mathbb{Q}$ then $(c_0,\cdots, c_{k-1}) \in \mathbb{Q}^k$.\\~\\
This completes the proof of lemma. $\square$

\section{Proof of theorem 1.2}
At first we prove theorem 1.2 and deduce theorem 1.1 as a corollary.\\
We shall use lemma 2.1 with $g(X) = P(X)$. Say, leading co-efficient of $P$ is $a_0 > 0$. \\
Let $k$ be an integer $\geq 2$.\\
Let $(c_0,\cdots,c_{k-1}) \in \mathbb{R}^{k}$ be the tuple as in lemma 2.1.\\
We continue to use notations from lemma 2.1. \\
Put $f(X)=F (X,c_0,\cdots,c_{k-2},c)= c_0X^{k-1}+\cdots+c_{k-2}X+c$ where $c$ is any real number satisfying $c < c_{k-1} < c+\frac{k}{k-1}$.\\
Now 
\[
\frac {1} {f(X)} -\frac {1} {f(X+1)} -\frac {1} {g(X+1)} =\frac {g(X+1)(f(X+1)- f(X)) - f(X)f(X+1)} {f(X)f(X+1)g(X+1)} \]
Consider the expression on the numerator.\\
By choice of $(c_0,..,c_{k-2})$ the coefficients of $X^{2k-2},\cdots,X^{k}$ vanishes (i.e. $p_i=q_i$ holds for all $0 \leq i \leq k-2$).\\
Coefficient of $X^{k-1}$ is $(q_{k-1}(c_0,\cdots,c_{k-2},c) - p_{k-1}(c_0,\cdots,c))$.\\
From proof of lemma 2.1 we have $q_{k-1}$ does not depends on $x_{k-1}$.\\
Hence $q_{k-1}(c_0,..,c_{k-2},c)=q_{k-1}(c_0,\cdots,c_{k-2},c_{k-1})=p_{k-1}(c_0,\cdots,c_{k-2},c_{k-1})$.\\
(By choice of the tuple$ (c_0,\cdots,c_{k-1}))$\\
Again from the proof of lemma 2.1 the coefficient of $x_{k-1}$ in $p_{k-1}$ is $2x_0$. Thus\\~\\
$\begin{aligned}
&q_{k-1}(c_0,\cdots,c_{k-2},c) - p_{k-1}(c_0,\cdots,c)\\ 
&= p_{k-1}(c_0,\cdots,c_{k-2},c_{k-1})-p_{k-1}(c_0,\cdots,c_{k-2},c) \\
&= 2c_0(c_{k-1}-c)
\end{aligned}$\\~\\
Hence the coefficient of the leading term of the polynomial in the numerator is $2c_0(c_{k-1}-c)$. But $c_0 = a_0(k-1) > 0$. So $2c_0(c_{k-1} -c) > 0$.\\
The coefficient of the leading term of the polynomial in the denominator is $c_0^2 > 0$.\\
Hence the is a large enough natural number $N_1$ such that for all $ i \geq N_1 $
\[\frac {1} {f(i+i_0)} -\frac {1} {f(i+i_0+1)} -\frac {1} {g(i+i_0+1)} > 0\]
i.e. \[\frac {1} {f(i+i_0)} -\frac {1} {f(i+i_0+1)} > \frac {1} {g(i+i_0+1)} \]\\
From here telescoping we get \\
\[\frac {1} {f(n+i_0)} > \sum_{i=n+1}^{\infty} \frac {1} {g(i+i_0)}\tag{3.1}\]
for all $n \geq N_1+i_0$. 
Now \\~\\
$\begin{aligned}&\frac {1} {f(X)+1} -\frac {1} {f(X+1)+1} -\frac {1} {g(X+1)} \\= &\frac {g(X+1)(f(X+1)- f(X)) - (f(X)+1)(f(X+1)+1)} {(f(X)+1)(f(X+1)+1)g(X+1)}\end{aligned}$\\
\vspace{.5cm}\\
Note that $(f(X)+1)(f(X+1)+1) - f(X)f(X+1)$ has degree equal to $(k-1)$. Hence co-efficient of $X^{k-1}$ in $(f(X)+1)(f(X+1)+1)$ is $p_{k-1}(c_0, \cdots c_{k-2}, c) + 2a_0$. Similar calculation suggests that the coefficient of the leading term in the numerator is $2c_0(c_{k-1} - c-1) - 2a_0 = 2c_0(c_{k-1} - c - \frac {k} {k-1}) < 0$.\\
But the coefficient of the leading term in the denominator is $c_0^2 > 0$.\\
Thus there is a large enough integer $N_2$ such that for all $i \geq N_2$
\[\frac {1} {f(i+i_0)+1} -\frac {1} {f(i+i_0+1)+1}  < \frac {1} {g(i+i_0+1)} \]
Again by telescoping \\
\[\frac {1} {f(n)+1} < \sum_{i=n+1}^{\infty} \frac{1}{g(i+i_0)}\tag{3.2}\]
for all $n \geq N_2 +i_0$. \\~\\
Put $N_0 = \max \{N_1, N_2\} + i_0$.\\
Then 
\[\frac {1} {f(n)+1} < \sum_{i=n+1}^{\infty} \frac {1} {g(i)} < \frac {1} {f(n)}\tag{3.3}\]
for all $n \geq N_0$. \\~\\
From $(3.3)$ the theorem 1.2 follows.\hfill $\square$\\~\\
\textbf{Remark 3.1: }i) The proof of lemma 2.1 and proof of theorem 1.1 gives an algorithm to compute the polynomial $f(X)$ mentioned in the statement of the theorem. Note that this polynomial is not unique.  We can choose infinitely many distinct values for the constant term. The integer $N_0$ depends on the choice of the polynomial $f$.\\
ii) It is natural to ask whether we can take $c=c_{k-1}$.\\
Put
\[f_1(X)=c_0X^{k-1}+c_1X^{k-2}+\cdots+c_{k-2}X+c_{k-1}.\]
We consider three possible cases:\\~\\
\textbf{Case I: } $p_i(c_0,\cdots,c_{k-1})=q_i(c_0,\cdots,c_{k-1})$ for all $0 \leq i \leq 2k-2$.\\
Then 
\[\frac {1} {g(X+1)}=\frac {1} {f_1(X)} - \frac {1} {f_1(X+1)}.\]\\
Then by telescoping $\sum_{i=n+1}^{\infty} \frac {1} {g(i+i_0)}=\frac {1} {f_1(n+i_0)}$ holds for all positive integer $n$.\\
Since $c_{k-1} < \frac{k} {k-1} + c_{k-1}$ arguments as in proof of the theorem implies $\sum_{i=n+1}^{\infty} \frac {1} {g(n+i_0)} > \frac {1} {f_1(n+i_0)+1}$ for large enough $n$.\\
Thus 
\[ f_1(n+i_0) \leq \Big(\sum_{i \geq n+1}\frac {1} {g(i+i_0)}\Big)^{-1} < f_1(n+i_0)+1 \] holds for large enough $n$.\\~\\
Now say, case I does not hold.
Then there is a $i$ with $0 \leq i \leq 2k-2$ such that $p_i \neq q_i$ at the point $(c_0,\cdots,c_{k-1})$.\\
Let $i_0$ be the minimum of such $i$ s. Note that here we must have $i_0 \geq k$. Now\\~\\
\textbf{Case II: } $q_{i_0}(c_0,\cdots,c_{k-1})  > p_{i_0}(c_0,\cdots,c_{k-1})$\\
Then $\sum_{i=n+1}^{\infty} \frac {1} {g(i+i_0)} < \frac {1} {f_1(n+i_0)}$ for large enough $n$ by arguments similar to the proof of the theorem.\\
Since $c_{k-1} < c_{k-1} + \frac {k} {k-1}$, so we have $\sum_{i=n+1}^{\infty} \frac{1} {g(n+i_0)} > \frac {1} {f_1(n+i_0)+1}$ for large enough $n$.\\
Thus $f_1(X)$ satisfies the required property of $f$ in the theorem.\\~\\
\textbf{Case III: }  $q_{i_0}(c_0,\cdots,c_{k-1})  < p_{i_0}(c_0,\cdots,c_{k-1})$\\
Again by similar arguments $\sum_{i=n+1}^{\infty} \frac {1} {g(i+i_0)} > \frac {1} {f_1(n+i_0)}$ for large enough $n$.\\
Now since $c_{k-1}-1 < c_{k-1}$, so we have $\sum_{i=n+1}^{\infty} \frac {1} {g(i+i_0)} < \frac {1} {f_1(n+i_0)-1}$ for large enough $n$. So we can not take $c = c_{k-1}$ but can take $c = c_{k-1}-1$.\\~\\
\textbf{Corollary 3.2 (Theorem 1.1):} Put $P(X) = X^k$, $i_0 = 0$. Using lemma-2.1 we conclude that $(c_0, \cdots, c_{k-1}) \in \mathbb{Q}^k$. Clearly one can choose $c$ to be rational. Hence theorem 1.1. Note that considerations in remark $3.1$ hold accordingly. 
\section{Computation of $a_{n}$}
First we make a small observation:\\~\\
\textbf{Remark 4.1 : } Let $f$ be any polynomial given by theorem 1.2. Let $N_0$ be the corresponding integer. Then from the inequality in theorem 1.2 it follows that for all $n \geq N_0$\\
i) $a_n$ is either $[f(n)]$ or $[f(n)] + 1$.\\
ii) If $f(n)$ is an integer for some $n$ then $a_n=f(n).$\\
iii) Conclusion in i) and ii) continue to hold if we replace the `$<$' sign in the inequality at the left hand side  in the statement of the theorem by `$\leq$'.\\~\\
For the rest of the section we shall assume that $P$ is defined over $\mathbb{Q}$. Further, we shift the polynomial so that we can take $i_0 = 0$. The polynomial $X^k$ satisfies these properties.
\subsection{A general algorithm }
Fix a polynomial $P$ as above.\\
Calculate $(c_0,\cdots,c_{k-1})$.\\
Write $c_i=\frac {u_i} {v_i}$ where $u_i,v_i$ are integers with $v_i > 0$ and $\gcd(u_i,v_i)=1$ for all $0 \leq i \leq k-1$.\\
Put $V= \text{lcm} \hspace{1mm} (v_0,\cdots,v_{k-2})$.\\
Now consider two cases:\\~\\
\textbf{Case I :} $v_{k-1}$ does not divide $V$.\\
Write $c_{k-1} = [c_{k-1}] + \frac {r_{k-1}} {v_{k-1}}$ where $r_{k-1}$ is a positive integer. Since $\gcd(u_{k-1}, v_{k-1}) = 1$, one has $\gcd(r_{k-1}, v_{k-1}) = 1.$ So $\frac {r_{k-1}} {v_{k-1}} \neq \frac {n} {V}$ for any $n \in \mathbb{Z}$.\\
Let $r \in \{0,\cdots,V-1\}$ be fixed.\\
Then there is a unique integer $n(r)$ such that $n(r) -\frac {r} {V} < c_{k-1} < n(r)+1 -\frac {r} {V}$.\\
Put
\[h(X)=c_0X^{k-1}+\cdots+c_{k-2}X=\frac {h_0(X)} {V}\]
where $h_0(X) \in \mathbb{Z}[X]$.\\
Let
\[f_r(X)=h(X)+n(r)-\frac {r} {V}.\]
Now there is an integer $N(r)$ such that 
\[f_r(n) < \frac {1} {\sum_{i=n+1}^{\infty} \frac {1} {P(i)}} < f_{r}(n)+1\] for all $n \geq N(r)$. Choose such a $N(r)$.\\
We do this for each $r \in \{0,\cdots,V-1\}$.\\
Put $N=\max \hspace{1mm} \{N(0),\cdots,N(V-1)\}.$\\
Let $n \geq N$ and $r$ be such that $r \in \{0,\cdots,V-1\}$ and $h_0(n)\equiv r\;(\text{mod} V)$. Clearly such $r$ exists and is unique.\\
Then using remark 4.1 $(ii)$ we have $a_n=f_r(n)$.\\
Now note that $n_1\equiv n_2 \,\text{mod}(V)$ implies $h_0(n_1)\equiv h_0(n_2)\;\text{mod}(V)$.\\
Thus in this case we have a closed form formula for $a_n$ depending on equivalence class of $n$ modulo $V$ whenever $n \geq N$ .\\~\\
\textbf{Case II: } Case I does not hold.\\
Fix $r \in \{1,\cdots,V\}$.\\
If $\frac{r} {V} \neq 1+[c_{k-1}] - c_{k-1}$, there is an unique integer $n(r)$ such that $n(r) - \frac {r} {V} < c_{k-1} < n(r) +1 - \frac {r} {V}$. \\
Otherwise there is an integer $n$ such that $n - \frac {r} {V}=c_{k-1}$.\\
Now we need to do further calculation and find out which case of remark 3.1 ii) holds.
If case I or II holds then put $n(r)=n$. If case III holds put $n(r)=n-1$.\\
Let
\[f_r(X)=h(X)+n(r)- \frac {r} {V}.\]
Using previous arguments and the discussion in remark 3.1 ii), there is an integer $N(r)$ such that
$$f_r(n) \leq \frac {1} {\sum_{i=n+1}^{\infty} \frac {1} {P(i)}} < f_{r}(n)+1$$ for all $n \geq N(r)$.\\ 
Due to remark 3.1 iii) the arguments of case I goes through from here.\\

\subsection{Explicit formulae}
Consider the polynomial $P(X) = X^k$, $k \geq 2$.\\~\\
For $k = 2$, $(c_0, c_1) = (1, \frac {1} {2})$.\\
Here $a_n = n$ for all $n \geq 1$ (ie one can take $N_0 = 1$).\\~\\
For $k = 3$, $(c_0, c_1, c_2) = (2, 2, 1)$.\\
Here $a_n = 2n(n+1)$ for all $n \geq 1$.\\~\\
For $k = 4$, $(c_0, c_1, c_2, c_3) = (3, \frac {9} {2}, \frac {15} {4}, \frac {9} {8}).$\\
Here 
$$ a_n = 
\begin{cases}
3X^3 + \frac {9} {2}X^2 + \frac {15} {4}X + 1 \; \text{if} \; n \equiv 0 (\text{mod}\, 4),\\
3X^3 + \frac {9} {2}X^2 + \frac {15} {4}X + \frac {3} {4} \; \text{if} \; n \equiv 1 (\text{mod}\, 4),\\
3X^3 + \frac {9} {2}X^2 + \frac {15} {4}X + \frac {1} {2} \; \text{if} \; n \equiv 2 (\text{mod}\, 4),\\
3X^3 + \frac {9} {2}X^2 + \frac {15} {4}X + \frac {1} {4} \; \text{if} \; n \equiv 3 (\text{mod}\, 4).
\end{cases}
$$\\
For $k = 5$, $(c_0, c_1, c_2, c_3, c_4) = (4, 8 , \frac {28} {3}, \frac {16} {3}, -\frac {2} {9})$.\\
Here 
\[a_n =
\begin{cases}
4X^4 + 8X^3 + \frac {28} {3}X^2 + \frac {16} {3}X - 1 \; \text{if} \; n \equiv 0 (\text{mod}\,3),\\
4X^4 + 8X^3 + \frac {28} {3}X^2 + \frac {16} {3}X - \frac {2} {3} \; \text{if} \; n \equiv 1 (\text{mod}\,3),\\
4X^4 + 8X^3 + \frac {28} {3}X^2 + \frac {16} {3}X - 1 \; \text{if} \; n \equiv 0 (\text{mod}\,3).
\end{cases}
\]\\
\textbf{Remark 3.2 :}  Results above answer two questions due to Kotesovec \cite{oeis}.

\section{Concluding remarks}
We end with some questions associated to the system of equations which come up in lemma 2.1 : \\
i) Consider the sequence of polynomials $\{P_{k}(X)\}_{k \geq 2}$ given by $P_{k}(X) = X^{k}$. With this sequence one can associate a sequence $(c_0, c_1, \cdots )$ where $c_{i}$ is a function $\mathbb{N} - \{1, \cdots i+1 \} \to \mathbb{R}$ such that $(c_0(k), \cdots, c_{k-1}(k))$ is the tuple associated to $P_k(X)$. From lemma 2.1 it follows that $c_{i}(k)$ must be a rational fuction of $k$. Computing first few elements of $(c_{0}(k), c_{1}(k), \cdots )$ one sees that it is actually a polynomial in $k$. This leads to the question if $c_i$ is always a polynomial in $k$.\\
ii) One can consider a sequence qiven by $P_{k}(X) = X^kP_0(X)$ for some fixed polynomial $P_0(X)$ and ask similar question.\\
iii) Fix two polynomials $P(X), Q(X)$. Construct a sequence by $P_{k}(X) = P(X)Q(X)^k$. In this case $c_i$ need not be a rational function but one may like to study the behaviour of the associated sequence of functions.\\  

\section{Acknowledgements}
This article is a refined version of a manuscript submitted to Journal of Integer Sequence (June, 2017). I am thankful to Prof. J. Shallit for some advice. This is part of a work which was developed over a long period of time (2014-2017) during my stay in Chennai Mathematical Institute (CMI). I am thankful to all my friends in CMI.

\end{document}